\documentstyle[12pt,oneside]{amsart}

\title{Jet schemes, arc spaces and the Nash problem}
\author{Shihoko Ishii} 
\address{Department of Mathematics, Tokyo Institute of
Technology, Oh-Okayama, Meguro, Tokyo, Japan
\newline
e-mail :  ishii.s.ac@@m.titech.ac.jp}

\newcommand{\bC}{{\Bbb C}}
\newcommand{\bP}{{\Bbb P}}

\newcommand{\bQ}{{\Bbb Q}}

\newcommand{\bR}{{\Bbb R}}
\newcommand{\bN}{{\Bbb N}}
\newcommand{\bA}{{\Bbb A}}

\newcommand{\Spec}{\operatorname{Spec}}

\newcommand{\Hom}{\operatorname{Hom}}

\newcommand{\st}{{\Spec k[[t]]}}
\newcommand{\stm}{{\Spec k[t]/(t^{m+1})}}
\newcommand{\sTm}{{\Spec K[t]/(t^{m+1})}}
\newcommand{\sT}{{\Spec K[[t]]}}
\newcommand{\tm}{{k[t]/(t^{m+1})}}

\let \cedilla =\c
\renewcommand{\c}[0]{{\mathbb C}}  

\renewcommand{\o}[0]{{\mathcal O}} 
\newcommand{\spec}[0]{\operatorname{Spec}}
\newcommand{\sing}[0]{\operatorname{Sing}}

\newcommand{\GL}{\operatorname{GL}}

\def\to {\longrightarrow}

\newtheorem{thm}{Theorem}[section]

\newtheorem{lem}[thm]{Lemma}

\newtheorem{prop}[thm]{Proposition}
\newtheorem{problem}[thm]{Problem}

\theoremstyle{definition}
\newtheorem{defn}[thm]{Definition}

\newtheorem{say}[thm]{}
\newtheorem{exmp}[thm]{Example}

\newtheorem{rem}[thm]{Remark}

\theoremstyle{remark}
\begin{document}
\maketitle
\footnote{partially supported by Grant-In-Aid of Ministry of Science 
and Education in Japan}
\footnote{AMS Subject Classification 2000: primary 14D15, secondary 
14D22, 14J17}

\begin{abstract}
  This paper is an introduction to the jet schemes and the  arc space 
  of an algebraic variety.
  We also introduce the Nash problem on arc families.
\vskip.5truecm
\noindent
R\'esum\'e. Ce papier constitue une introduction aux espaces de jets
et \`a l'espace d'arcs d'une vari\'et\'e alg\'ebrique.
Nous introduisons \'egalement le probl\`eme de Nash 
pour les familles d'arcs.
\vskip.5truecm

\noindent
Keywords: arc space, jet scheme,  Nash problem
\end{abstract}

%%%%%%%%%%%%%%%%%%%%%%%%%%%%%%%%%%%%%%%%%%%%%%%%%f
\section{Introduction}

\noindent
  The concepts  jet scheme and  arc space  over an algebraic variety 
 or an analytic space is introduced by Nash 
 in his preprint in 1968 
 which is later  published as \cite{nash}.
  The study of these spaces was further developed by 
   Kontsevich,  Denef and  Loeser 
 as the theory 
  of motivic integration, see \cite{ko}, \cite{ DL1}, \cite{ DL2}, 
  \cite{ DL3}, \cite{DL4},\cite{ DL5}. 
   These spaces  
  are considered as something to represent
  the nature of the singularities of the base space. 
 In fact, papers \cite{ein}, \cite{e-Mus}, \cite{must01}, \cite{must02} 
  by  Musta\cedilla{t}\v{a}, Ein and Yasuda show that  geometric 
  properties of the jet schemes determine certain properties 
  of the singularities 
  of the base space. 

  In this paper, we provide the beginners with the basic knowledge  of 
  these spaces and the Nash problem.
  One of   powerful  tools to work on these space is the motivic 
  integration.
  But this paper does not step into this theory, 
  as there are already very good introductory papers on 
  the motivic integration by A. Craw \cite{craw}, W. Veys 
  \cite{veys} and F. Loeser \cite{l}.
    We devote into the  basic study of  geometric structure of arc 
    spaces and jet schemes.  
  We also give an introduction to the Nash problem which was posed 
  in \cite{nash}.   

  Throughout this paper the base field \( k \) is algebraically closed 
  field of arbitrary characteristic and a variety is an irreducible 
  reduced scheme of finite type over \( k \).
  A scheme of finite type over \( k \) is always separated over \( k \).

  We omit the proofs of  statements whose references are thought 
  to be easily accessible. 
  We assume the reader to have knowledge in the Hartshorne's textbook 
  \cite{ha}.

  The author expresses  her hearty thanks to Clemens Bruschek who
  read the preliminary version of this paper and asked many constructible 
  questions.
  In order to answer his questions, many parts were improved.

%%%%%%%%%%%%%%%%%%%
\section{Construction of jet schemes and arc spaces}

\begin{defn}
  Let \( X \) be a scheme of finite type over \( k \)
and $K\supset k$ a field extension.
  For \( m\in \bN \), a  \( k \)-morphism \( \Spec K[t]/(t^{m+1})\to X \) is called an \( m \)-jet
  of \( X \) and a  \( k \)-morphism \( \Spec K[[t]]\to X \) is called an 
  {\it arc} of \( X \).
  We denote the unique point of \( \sTm \) by \( 0 \), 
  while the closed point of \( \Spec K[[t]] \) by \( 0 \) and
  the generic point by \( \eta \).
\end{defn}

\begin{thm}
\label{existence}
  Let \( X \) be a scheme of finite type over \( k \).
  Let \( {\cal S}ch/k \) be the category of \( k \)-schemes  
   and \( {\cal S}et \) the category of sets.
  Define a contravariant functor  \( F^X_{m}: {\cal S}ch/k \to {\cal S}et \)
  by 
$$
 F^X_{m}(Z)=\Hom _{k}(Z\times_{\Spec k}\Spec k[t]/(t^{m+1}), X).
$$
  Then, \( F^X_{m} \) is representable by a scheme \( X_{m} \) of finite
  type over \( k \), that is
$$
 \Hom _{k}(Z, X_{m})\simeq\Hom _{k}(Z\times_{\Spec k}
\Spec k[t]/(t^{m+1}), X).
$$ 
   This \( X_{m} \) is called the {\it space of \( m \)-jets} of \( X \) 
   or the {\it \( m \)-jet scheme} of \( X \).
\end{thm}

  This proposition is proved in  \cite[p. 276]{blr}.
  In this paper, we prove this by a concrete construction for  affine 
  \( X \)
  first and then patching them together for a general \( X \).
  For our proof, we need some preparatory discussions.
  
\begin{say}    
  Let \( X \) be a \( k \)-scheme. 
  Assume that \( F_{m}^X \) is representable by \( X_{m} \) for every
   \( m\in \bN \).
  Then, for \( m<m' \), the canonical surjection \( k[t]/(t^{m'+1})
  \to k[t]/(t^{m+1}) \)
  induces a morphism \[ \psi_{m',m}:X_{m'}\to X_{m}. \]  
  Indeed, the canonical surjection 
  \( k[t]/(t^{m'+1})\to k[t]/(t^{m+1}) \)
  induces a morphism 
  \[ Z\times_{\spec k}\spec k[t]/(t^{m'+1})\leftarrow Z\times_{\spec 
  k}\stm , \]
  for an arbitrary \( k \)-scheme \( Z \).
  Therefore we have a map \( \Hom_{k}(Z\times_{\spec k}\spec k[t]/(t^{m'+1})
  , X) \to \Hom_{k}(Z\times_{\spec   k}\stm , X) \) which gives the map
  \[ \Hom_{k}(Z, X_{m'})\to \Hom_{k}(Z, X_{m}). \]
  Take, in particular, \( X_{m'} \) as \( Z \), 
  \[ \Hom_{k}(X_{m'},X_{m'})\to \Hom_{k}(X_{m'},X_{m}) \]
  then the image of \( 
  id_{X_{m'}}\in \Hom (X_{m'},X_{m'}) \) 
  by this map gives the required morphism.

  This morphism \( \psi_{m',m} \) is called a {\it truncation map}.
  In particular for \( m=0 \), 
  \( \psi_{m',0}:X_{m'}\to X  \)
  is denoted by \( \pi_{m} \).
  When we need to specify the scheme \( X \), 
  we denote it by \( \pi_{Xm} \).
  
  Actually \( \psi_{m',m} \) ``truncates'' a power series in the following 
  sense:
  A point \( \alpha \)  of \( X_{m'} \)  gives an \( m' \)-jet \( \alpha:
  \Spec K[t]/(t^{m'+1})
  \to X \), which corresponds to a ring homomorphism \( \alpha^*: A\to 
  K[t]/(t^{m'+1}) \), where \( A \) is the affine coordinate ring of 
  an affine neighborhood of the image of \( \alpha \).
  For every \( f\in A \), let \[ 
  \alpha^*(f)=a_{0}+a_{1}t+a_{2}t^2+\cdots +a_{m}t^m+ \cdots 
  +a_{m'}t^{m'}, 
  \]
  then 
  \[ (\psi_{m',m}¥(\alpha))^*(f)= a_{0}+a_{1}t+a_{2}t^2+\cdots +a_{m}t^m.\]
  This fact can be seen by letting \( Z=\{\alpha\} \) in the above 
  discussion.

  As we did already in the above argument, 
  we denote the point of \( X_{m} \) corresponding to \( \alpha:
  \Spec K[t]/(t^{m+1})
  \to X \)
  by the same symbol \( \alpha \).
  Then, we should note that \( \pi_{m}(\alpha)=\alpha(0) \).
\end{say}

\begin{prop}
  Let \( f:X\to Y \) be a morphism of \( k \)-schemes of finite type.
  Assume that the functor \( F_{m}^X \) and \( F_{m}^Y \) are 
  representable by \( X_{m}  \) and \( Y_{m} \), respectively.
  Then the canonical morphism \( f_{m}:X_{m}\to Y_{m} \) 
  is induced for every \( m\in \bN \) such that the 
  following diagram is commutative:
  \[ \begin{array}{ccc}
      X_{m}& \stackrel{f_{m}}\longrightarrow & Y_{m}\\
      \pi_{Xm} \downarrow\ \ \ \ \ & & \ \ \ \ \downarrow \pi_{Ym}\\
      X & \stackrel{f}\longrightarrow & Y\\
      \end{array}. \]
\end{prop}

\begin{pf}
    Let \( X_{m}\times \stm\to X \) be the ``universal family'' of \( m 
  \)-jets of \( X \), i.e., it corresponds to the identity map in 
  \( \Hom_{k}(X_{m},X_{m}) \).
  By compositing this map and \( f:X\to Y \), we obtain a morphism
  \[  X_{m}\times \stm\to Y, \]
  which gives a morphism \( X_{m}\to Y_{m} \).
  Pointwise, this morphism maps  an \( m \)-jet \( \alpha\in X_{m} \) 
  of \( X \) to
  the composite \( f\circ \alpha \) which is an \( m \)-jet  of \( Y \).
  To see this, just take a point \( \alpha\in X_{m} \) and see the 
  image of \( \{\alpha\}\times \stm \to Y \).
  The commutativity of the diagram follows from this description.
\end{pf}

\begin{prop}
\label{etale}
%%**(etale)
  For  \( k \)-schemes \( X \) and \( Y \),
  assume that the functor \( F_{m}^X \) and \( F_{m}^Y \) are 
  representable by \( X_{m}  \) and \( Y_{m} \), respectively.
 If \( f:X\to Y \) is an \'etale morphism, then \( X_{m}\simeq 
 Y_{m}\times_{Y}X \), for every \(m\in \bN \).
\end{prop}

\begin{pf}
  By the above proposition we have a commutative diagram:
  \[ \begin{array}{ccc}
  X_{m} & \stackrel{f_{m}}\longrightarrow & Y_{m}\\
  \downarrow & & \downarrow \\
  X & \stackrel{f}\longrightarrow & Y\\
  \end{array}. \]
  It is sufficient to prove that for every commutative diagram:
    \[ \begin{array}{ccc}
  Z & \longrightarrow & Y_{m}\\
  \downarrow & & \downarrow \\
  X & \stackrel{f}\longrightarrow & Y\\
  \end{array}, \]
  there is a unique morphism  \( Z\to X_{m} \) which is compatible 
  with the projections to \( X \) and \( Y_{m} \).
  Now we are given the following commutative diagram:
  \[ \begin{array}{ccc}
  Z & \longrightarrow & Z\times_{\spec k}\spec k[t]/(t^{m+1})\\
  \downarrow & & \downarrow \\
  X & \stackrel{f}\longrightarrow & Y\\
  \end{array}. \]
  As \( f \) is \'etale, there is a unique morphism 
  \( Z\times_{\spec k}\spec k[t]/(t^{m+1})\to X \) which makes the two
  triangles commutative.  
  This gives the required morphism:

%   As \( \varprojlim_{m}(Y_{m}\times_{Y}X)=(\varprojlim_{m}Y_{m})\times_{Y}X 
%   \), it is sufficient to prove the assertion for \( m\in \bN \).
%   By the commutative diagram:
%   \[ \begin{array}{ccc}
%   X_{m} & \stackrel{f_{m}}\longrightarrow & Y_{m}\\
%   \downarrow & & \downarrow \\
%   X & \stackrel{f}\longrightarrow & Y\\
%   \end{array}, \]
%   we obtain a morphism \(\varphi: X_{m}\to Y_{m}\times_{Y}X \).
%   On the other hand, the projection \( Y_{m}\times_{Y}X\to Y_{m}  \) 
%   corresponds to a morphism \(  Y_{m}\times_{Y}X \times 
%   _{\spec k }\stm \to Y \) which is a part of the following commutative
%   diagram:
%   \[ \begin{array}{ccc}
%   Y_{m}\times_{Y}X & \hookrightarrow &  Y_{m}\times_{Y}X \times 
%   _{\spec k }\stm\\
%   \downarrow & & \downarrow\\
%   X & \stackrel{f}\longrightarrow &Y\\
%   \end{array}. \]
%   As \( f \) is  \'etale, there is a morphism \\
%   \(   Y_{m}\times_{Y}X \times 
%   _{\spec k }\stm\to X \) which makes  the  two triangles commutative. 
%   The corresponding morphism \(  Y_{m}\times_{Y}X \to X_{m} \)
%   is the inverse morphism of \( \varphi \).  
\end{pf}

As a corollary of this proposition, we obtain the following lemma:

\begin{lem}
\label{open}
%%*(open)
  Let  \( U\subset X \) be an open subset of a \( k \)-scheme \( X \).
  Assume the functors \( F^X_{m} \) and \( F^U_{m} \) are representable 
  by \( X_{m} \) and \( U_{m} \), respectively.
  Then,
  \(U_{m}= \pi_{Xm}^{-1}(U) \).
\end{lem}  
\vskip.5truecm
[{\it Proof of Theorem \ref{existence}}]
  Since a \( k \)-scheme \( X \) is separated, the intersection of 
  two affine open subsets is again affine.
  Therefore, by Lemma \ref{open}, it is sufficient to prove the 
  representability of \( F_{m}^X \) for affine \( X \).  
  Let \( X \) be \( \spec R \), where we denote 
  \( R=k[x_{1},\ldots,x_{n}]/(f_{1},..,f_{r}) \).
  It is sufficient to prove the representability for an affine variety
   \( 
  Z=\spec A \).
  Then, we obtain that
  \[(\ref{existence}.1)\ \ \ \ \ \Hom (Z\times \spec k[t]/(t^{m+1}), X)\simeq
  \Hom (R, A[t]/(t^{m+1})) \]
  \[\simeq \{\varphi\in \Hom(k[x_{1},.,x_{n}], A[t]/(t^{m+1})) \mid
  \varphi(f_{i})=0 \ \operatorname{for} \ i=1,..,r\}. \]
  If we write \( 
   \varphi(x_{j})=
   a_{j}^{(0)}+a_{j}^{(1)}t+a_{j}^{(2)}t^2+..+a_{j}^{(m)}t^m \) for 
   \( a_{j}^{(l)}\in A \),
   it follows that  \[ \varphi(f_{i})= 
   F_{i}^{(0)}(a_{j}^{(l)})+ F_{i}^{(1)}(a_{j}^{(l)})t +..+ 
   F_{i}^{(m)}(a_{j}^{(l)})t^m\]
   for  polynomials  \(  F_{i}^{(s)} \) in \( a_{j}^{(l)} \)'s. 
  Then the above set (\ref{existence}.1)  is represented as follows:
  \[ =\{\varphi\in \Hom (k[x_{j},x_{j}^{(1)},.,x_{j}^{(m)}\mid 
  j=1,.,n], A) \mid \varphi(x_{j}^{(l)})=a_{j}^{(l)},\ 
  F_{i}^{(s)}(a_{j}^{(l)})=0\} \]
  \[ =\Hom (k[x_{j},x_{j}^{(1)},.,x_{j}^{(m)}]/( 
  F_{i}^{(s)}(x_{j}^{(l)})),
  A). \]    
  If we write 
  \( X_{m}=\spec k[x_{j},x_{j}^{(1)},.,x_{j}^{(m)}]/
  ( F_{i}^{(s)}(x_{j}^{(l)})) \),
  the last set  is bijective to
  \[ \Hom (Z, X_{m}). \]
\( \Box \)

\begin{rem}
The functor \( F_{m}^X \) is also representable even for \( k 
\)-scheme of non-finite type over \( k \).
The existence of jet schemes for wider class of schemes is presented in 
\cite{voj}.
%   The defining equations \( F_{i}^s(x_{j}^{(l)}) \)'s of \( X_{m} \)
%   are obtained as 
%   follows:
%   Let \( D \) be a derivation of \( k[x_{j},x_{j}^{(1)},.,x_{j}^{(m)}] \)
%   defined by \( D(x_{j}^{(l)})=x_{j}^{(l+1)} \), where we define 
%   \( x_{j}^{(l)}=0 \) for \( l>m \).
% %
%   Then,   it follows that \(  F_{i}^s(x_{j}^{(l)})=D^s(f_{i}) \).
\end{rem}

\begin{exmp}
  For \( X=\bA_{k}^n \), it follows \( X_{m}=\bA_{k}^{n(m+1)} \).
  Indeed, this is the case that all \( f_{i}=0 \), therefore
  all \( F_{i}^{(s)}=0  \) in the proof of 
 Proposition \ref{existence}.
\end{exmp}

\begin{exmp}
  Let \( X \) be a hypersurface in \( \bA_{k}^3 \) defined by 
  \( f=xy+z^2=0 \).
  Then, \( X_{2} \) is defined in \( \bA_{k}^{9} \) by 
  \( xy+z^2=x^{(1)}y+xy^{(1)}+2zz^{(1)}=x^{(2)}y+x^{(1)}y^{(1)}
  +xy^{(2)}
  +z^{(1)}z^{(1)}+2zz^{(2)}=0 \).
  One can see that \( X_{2} \) is irreducible and not normal.
  Indeed, as \( X\setminus \{0\} \) is non-singular, \( \pi_{2}^{-1}(X
  \setminus \{0\}) \) is a 6-dimensional irreducible variety. 
  On the other hand \( \pi_{2}^{-1}(0) \) is a hypersurface in \( \bA^6 
  \), and therefore it is of dimension 5. 
  Since \( X_{2} \) is defined by 3 equations, every irreducible 
  component of \( X_{2} \) has dimension \( \geq 9-3=6 \).
  By this, \( \pi_{2}^{-1}(0) \) does not produce an irreducible component 
  of \( X_{n} \), which yields the irreducibility of \( X_{m} \).
  Looking at the Jacobian matrix, one can see that the singular locus 
  of \( X_{2} \) is \( \pi_{2}^{-1}(0) \) which is of codimension one 
  in \( X_{2} \). 
  Therefore, \( X_{2} \) is not normal. 
\end{exmp}

  Let \( X_{1} \) be the 1-jet scheme of \( X \). 
  Then for every closed point \( x\in X \), the set of  closed points of 
  \( \pi_{1}^{-1}(x) \) is the set of 
  morphisms \( \spec k[t]/(t^2)\to X \) with the image \( x \).
  This set is the Zariski tangent space of \( X \) at \( x \).
  Therefore, we can regard \( X_{1} \) as the  ``tangent bundle'' 
  of \( X \).

\begin{exmp}
\label{nonirred}
%*(nonirred)  
  Let \( X \) be a curve defined by \( x^2-y^2-x^3=0 \) in \( \bA_{k}^2 \).
  Then \( \pi_{1}^{-1}(X\setminus\{0\})\to X\setminus \{0\} \) is \( \bA_{k}^1 
  \)-bundle, 
  therefore \( \pi_{1}^{-1}(X_{reg}) \) is 2-dimensional. 
  On the other hand \(\pi_{1}^{-1}(0) \simeq \bA_{k}^2 \).
  Hence, \( X_{1} \) has two irreducible components, 
  \( \overline{\pi_{1}^{-1}(X_{reg})} \) and \( \pi_{1}^{-1}(0) \).
 \end{exmp}
 
\begin{defn}
  The system \( \{\psi_{m', m}:X_{m'}\to X_{m}\}_{m< m'} \) is a 
  projective system.  
  Let \( X_{\infty}=\varprojlim _{m}X_{m} \) and call it the 
{\it space of 
  arcs} of \( X \) or {\it arc space} of \( X \).
Note that \( X_{\infty} \) is not of finite type over $k$
if \( \dim X>0 \).
\end{defn}

\begin{rem}
  One may be afraid that the projective limit scheme
  \( \varprojlim _{m}X_{m} \) may not exist.
  But in our case we need not to worry, since for an affine scheme 
  \( X =\spec A\), the \( m \)-jet scheme \( X_{m}=\spec A_{m} \) 
  is affine for every 
  \( m \). 
  Here, the  morphisms \( \psi^*_{m',m}:A_{m}\to A_{m'} \) 
  corresponding to \( \psi_{m',m} \) are direct 
  system.
  It is well known that there is a direct limit \( 
  A_{\infty}=\varinjlim_{m}A_{m} \) in the category of \( k 
  \)-algebras.
  The affine scheme \( \spec A_{\infty} \) is our projective limit of 
  \( X_{m} \).
  For a general \( k \)-scheme \( X \), we have only to patch affine 
  pieces \( \spec A_{\infty} \).
\end{rem}

  Using  the representability of \( F_{m}^X \) 
we obtain the following universal property of $X_{\infty}$:

\begin{prop}
\label{ft}
%*(ft)
  Let \( X \) be a scheme of finite type over \( k \).
  Then
  \[ \Hom _{k}(Z, X_{\infty})\simeq\Hom _{k}(Z\widehat\times_{\Spec k}\st, X) \]
   for an arbitrary \( k \)-scheme \( Z \),
   where \( Z\widehat\times_{\Spec k}\st \) means the formal completion 
   of \( Z\times_{\Spec k}\st \) along the subscheme 
   \( Z\times _{\Spec k} \{0\} \).
\end{prop}

\begin{pf}
  By the representability of \( F_{m}^X \) we obtain an isomorphism of projective 
  systems:
  
  \[ \begin{array}{ccc}
  \downarrow& & \downarrow\\
 \Hom _{k}(Z, X_{m+1})&\simeq&\Hom _{k}(Z\times_{\Spec k}
\Spec k[t]/(t^{m+2}), X)\\
\downarrow& & \downarrow\\
 \Hom _{k}(Z, X_{m})&\simeq&\Hom _{k}(Z\times_{\Spec k}
\Spec k[t]/(t^{m+1}), X)\\
\end{array}
 .\]
 Then, we obtain an isomorphism of the projective limits:
 \[ \Hom_{k}(Z, \varprojlim _{m}X_{m})\simeq \Hom_{k} 
 (\varinjlim_{m}(Z\times_{\spec k}\spec k[t]/(t^{m+1})), X), \]
 which gives the required isomorphism.  
\end{pf}

\begin{rem}
  Consider the isomorphism of Proposition \ref{ft} in particular the 
  case \( Z=\spec A\) for a \( k \)-algebra \( A \), we obtain 
  \[ \Hom_{k}(\spec A, X_{\infty})\simeq \Hom_{k}(\spec A[[t]], X). \]
  Here, we note that in general
  \[ A\otimes_{k}k[[t]]\not\simeq A[[t]]\simeq 
  A\widehat{\otimes}_k  k[[t]], \]
  where \( A\widehat{\otimes}_k  k[[t]] \) is the completion of \( 
  A\otimes_{k}k[t] \) by the ideal \( (t) \).  
  Indeed, for example, for \( A=k[x] \), the ring
  \( A[[t]] \) contains  \( \sum_{i=0}^\infty f_{i}(x)t^i \) 
  such that  \( \deg f_{i} \) are unbounded,
  while   \(  A\otimes_{k}k[[t]] \) does not contain such an element.
  
  Now, consider the case \( A=K \) for an extension field \( K\supset k \),
  the bijection 
   \[ \Hom_{k}(\spec K, X_{\infty})\simeq \Hom_{k}(\spec K[[t]], X) \]
   shows that a \( K \)-valued point of \( X_{\infty} \) is an arc
   \( \sT\to X \). 
\end{rem}

\begin{defn}
  Denote the canonical projection \( X_{\infty }\to X_{m} \) induced 
  from the surjection \( k[[t]]\to k[t]/(t^{m+1}) \) by \( 
  \psi_{m} \) and the composite \( \pi_{m}\circ \psi_{m} \) by \( \pi \).
  When we need to specify the base space \( X \), we write it by \( 
  \pi_{X} \).
\end{defn}
  
    A point \( x \in X_{\infty} \)  
   gives an arc \( \alpha_{x}:\sT
  \to X \) and \( \pi(x)=\alpha_{x}(0) \), 
  where \( K \) is the residue field at \( x \). 
  As the case of  \( m \)-jets, 
   we denote  both \( x \in X_{\infty}  \) and \( \alpha 
  \) corresponding to \( x \) by the same symbol \( \alpha \).
  
  For every \( m\in \bN \), \( \psi_{m}(X_{\infty}) \) is a 
  constructible set, 
  since \( \psi_{m}(X_{\infty})=\psi_{m',m}(X_{m'}) \) for 
  sufficiently big \( m' \) (\cite{gr}).

\begin{defn} 
  Denote the canonical morphism  \( X\to X_{m} \) induced from the 
  inclusion \( k\hookrightarrow k[t]/(t^{m+1}) \) \( (m\in \bN\cup 
  \{\infty\}) \) by \( \sigma_{m} \).
  Here, we define  \( k[t]/(t^{m+1})= k[[t]] \) for \( m=\infty \).
  As \(  k\hookrightarrow k[t]/(t^{m+1}) \) is a section of the 
  projection \(   k[t]/(t^{m+1}) \to k\), 
  our morphism \( \sigma_{m}:X\to X_{m} \) is a section of 
  \( \pi_{m}:X_{m}\to X \).

\end{defn}

  For a point \( x\in X \), let \( K \) be the residue field at \( 
  x \), then \( \sigma_{m}(x) : \sTm \to X \) is an \( m \)-jet which 
  factors through \( \spec K\to X \) whose image is \( x \).
  Therefore, \( \sigma_{m}(x) \) is the constant \( m \)-jet at \( x \), 
  this is denoted my \( x_{m} \).

\begin{exmp}
  If \( X=\bA_{k}^n \), then \( X_{\infty}=\spec 
  k[x_{j},x_{j}^{(1)},x_{j}^{(2)}\ldots\mid j=1,\ldots, n] \) which is 
  isomorphic to  
  \( \bA_{k}^{\infty}=\spec k[x_{1},x_{2},\ldots, x_{i},\ldots] \).
  Here, we note that the set of closed points of \( \bA_{k}^{\infty} 
  \) does not necessarily coincide with the set \[ k^{\infty}:=
  \{(a_{1},a_{2},\ldots)\mid a_{i}\in k\}\]
 (see the following theorem).
\end{exmp}

\begin{thm}[\cite{i}, Proposition 2.10, 2.11]
\label{countable}
%***(countable)
  Every closed point of \( \bA_{k}^{\infty} \) is a \( k \)-valued 
  point if and only if \(  k \) is not a countable field.
\end{thm}

  The concept ``thin'' in the following is first introduced in 
  \cite{ELM}.

\begin{defn}
   Let \( X \) be a variety over \( k\). 
   We say that an arc \( \alpha:\sT  \to  X \) is {\it thin} 
   if \( \alpha  \) factors through a proper closed subvariety of \( X \).
  An arc which is not thin is called a {\it fat arc}.

  An irreducible subset \( C \) in \( X_{\infty} \)  is called a {\it 
  thin set} if \( C \) is contained in \( Z_{\infty} \) for a proper 
  closed subvariety \( Z\subset X \).
  An irreducible subset  in \( X_{\infty} \) which is not thin  is called a {\it 
  fat set}.
  
  In case an irreducible subset \( C \) has the generic point \( \gamma\in C 
  \) (i.e., the closure \( \overline{\gamma}  \) contains \( C \)),
  \( C \) is a fat set if and only if \( \gamma \) is a fat arc.
\end{defn} 

The following is proved in \cite[Proposition 2.5]{i2}:

\begin{prop}[\cite{i2} Proposition 2.5]
\label{fatbasic}
%*(fatbasic)
  Let \( X \) be a variety over \( k\) and \( \alpha:\sT\to X \) an 
  arc. 
  Then, the following hold:
\begin{enumerate}
\item[(i)]
  \( \alpha \) is a fat arc if and only if
  the ring homomorphism \( \alpha^*: \o_{X, \alpha(0)}\to K[[t]] \) 
  induced from \( \alpha \) is injective;
\item[(ii)]
  Assume that \( \alpha \) is fat.
  For an arbitrary proper birational morphism  
  \( \varphi:Y\to X \),  \( \alpha \) is lifted to \( Y \).
\end{enumerate}  
\end{prop}

\begin{rem}
  A fat set in \( X_{\infty} \) for a variety \( X \) introduces a 
  discrete valuation on the rational function field \( K(X) \) of \( X \).
  We do not give the construction of the valuation here.
  The reader may refer \cite{i2}.
  A Nash component (see the next section) is a fat set and the Nash 
  map (see the next section) is just the correspondence to associate a fat set to the 
  valuation induced from the fat set (\cite{i2}).
\end{rem}

\begin{exmp}
\label{fatex}
%*(fatex)
  One of  typical examples of  fat sets is an irreducible {\it 
  cylinder} ( i.e., the pull back \( \psi_{m}^{-1}(S) \) of a 
  constructible set \( 
  S\subset X_{m} \)) for a non-singular \( X \).
  Actually, take an \( m \)-jet \( \alpha_{m}:\stm\to X \) in \( C \), 
  then, at a  
  neighborhood of \( x=\alpha_{m}(0)=\pi_{m}(\alpha_{m}) \), \( X \) 
  is \'etale over \( \bA_{k}^n \).
  Therefore, we may assume that \( X=\bA_{k}^n \) and \( x=0 \).
  Assume that \( \psi_{m}^{-1}(\alpha_{m}) \) is thin, 
  then it is contained in \( Z_{\infty} \) for some proper closed 
  subset \( Z\in X \).
  Let the \( m \)-jet \( \alpha_{m} \) corresponds to a ring homomorphism 
  \[ \alpha_{m}^*: k[x_{1},\ldots,x_{n}]\to k[t]/(t^{m+1}), \ 
  \alpha_{m}^*(x_{i})=\sum _{j=1}^m a_{i}^{(j)}t^j .\]
  Let \( x_{i}^{(j)} \) be an indeterminate for every \( i=1,\ldots,n  \)
  and \( j\geq m+1 \).
  Let \[ \alpha^*:k[x_{1},\ldots,x_{n}]\to k(x_{i}^{(j)}\mid i=1,..,n, 
  j\geq m+1)[[t]] \]
  be an arc defined by \[ 
  \alpha^*(x_{i})=\sum_{j=1}^ma_{i}^{(j)}t^j+\sum_{j=m+1}^\infty x_{i}^{(j)}t^j. \]
  Let \( 
  \alpha^*(f)=F_{0}(a_{i}^{(j)},x_{i}^{(j)})+F_{1}(a_{i}^{(j)},x_{i}^{(j)})t+\cdots +
  F_{\ell}(a_{i}^{(j)},x_{i}^{(j)})t^{\ell}+\cdots \) for \( f\in I_{Z} \).
  Then, as \( x_{i}^{(j)} \)'s are indeterminates there is \( \ell \) such 
  that \( F_{\ell}\neq 0 \). 
  Hence, we obtain \( \alpha\in \psi_{m}^{-1}(C)\) such that \( 
  \alpha\not\in Z_{\infty} \).
\end{exmp}

\begin{exmp}[\cite{DEI}]
  For a singular variety \( X \), an irreducible cylinder is not necessarily 
  fat.
  Indeed, let \( X \) be the Whitney Umbrella that is a hypersurface defined by 
   \( xy^2-z^2=0 \) in \( \bA_{k}^3 \).
  For \( m\geq 1 \),
  let \[ \alpha_{m}^*:k[x,y,z]/(xy^2-z^2)\to k[t]/(t^{m+1}) \]
  be the \( m \)-jet defined by \( \alpha_{m}(x)=t, \alpha_{m}(y)=0, 
  \alpha_{m}(z)=0 .\)
  Then, the cylinder \( \psi_m^{-1}(\alpha_{m}) \) is contained in \(
  \sing(X)_{\infty} \), where \( \sing(X)= (y=z=0) \).
  This is proved as follows:
  Let an arbitrary \( \alpha\in \psi_{m}^{-1}(\alpha_{m}) \) be 
  induced from 
  \[ \alpha^*:k[x,y,z]\to k[[t]] \] with
  \[ \alpha^*(x)=\sum _{j=1}^\infty a_{j}t^j ,
   \alpha^*(y)=\sum _{j=1}^\infty b_{j}t^j ,
  \alpha^*(z)=\sum _{j=1}^\infty c_{j}t^j,  \]
  where we note that \( a_{1}=1 \).
  Then, the condition \( \alpha^*(xy^2-z^2)=0 \) implies 
  that the initial term of 
  \( \alpha^*(xy^2) \) and that of \( \alpha^*(z^2) \) cancel each 
  other. 
  If \( \alpha^*(y)\neq 0 \), then the order of \( \alpha^*(xy^2) \) 
  is odd, while if \( \alpha^*(z)\neq 0 \), the order of 
  \( \alpha^*(z^2) \) is even.
  Hence if \( \alpha^*(y)\neq 0 \) or \( \alpha^*(z)\neq 0 \), then 
  the initial term of 
  \( \alpha^*(xy^2) \) and that of \( \alpha^*(z^2) \) do not cancel 
  each other. 
  Therefore, \( \alpha^*(y)=\alpha^*(z)=0 \), which shows that 
  \( \psi_{m}^{-1}(\alpha_{m})\subset \sing(X)_{\infty} \).
  
\end{exmp}

%%%%%%%%%%%%%%%%%%%%%%%%%%%%%%%%%%%%%%%%%%%%%%%%%%%%%%%%%%%%%%
%%%%%%%%%%%%%%%%%%%%%%%%%%%%%%%%%%%%%%%%%%%%%%%%%%%%%%%%%%%%%%
\section{Properties of jet schemes and arc spaces}

\begin{say}
\label{invariant}
%***(invariant)
  Consider \( {\Bbb{G}}_m=\bA^1\setminus \{0\}=\spec k[s,s^{-1}] \) as a multiplicative 
  group scheme. 
  For \( m\in \bN\cup\{\infty\} \), the morphism \( \tm \to 
  k[s,s^{-1}, t]/(t^{m+1}) \) defined by \( t\mapsto s\cdot t \) 
  gives an action  \[ \mu_{m}: {\Bbb{G}}_m \times_{\spec k} \stm\to \stm\] of 
  \( {\Bbb{G}}_m  \) on \( \stm \).
  Therefore, it gives an action 
  \[ \mu_{Xm}: {\Bbb{G}}_m\times_{\spec k}X_{m}\to X_{m} \]
  of \( {\Bbb{G}}_m \) on \( X_{m} \).
  As \( \mu_{m} \) is extended to a morphism:
  \( \overline{\mu}_{m}:\bA^1\times _{\spec k} \stm\to \stm  \), we 
  obtain the extension 
   \[ \overline{\mu}_{Xm}: \bA^1\times_{\spec k}X_{m}\to X_{m} \]
   of \( \mu_{Xm} \).
   
   Note that \( \overline{\mu}_{Xm}(\{0\}\times \alpha)=x_{m} \),
   where \( x_{m} \) is the trivial \( m \)-jet on \( x=\alpha(0)\in X \). 
   Therefore, every orbit \( \mu_{Xm}({\Bbb{G}}_m\times \{\alpha\}) \) contains the 
   trivial \( m \)-jet on \( \alpha(0) \) in its closure.
\end{say}

\begin{prop}
  For \( m\in \bN\cup\{\infty\} \), let  \( Z \subset X_{m}\) be  
  an \( {\Bbb{G}}_m \)-invariant 
  closed subset.  
  Then the image \( \pi_{m}(Z) \) is closed in \( X \). 
  In particular the image \( \pi_{m}(Z) \) of an irreducible component of 
  \( Z\subset  X_{m} \)  is closed in \( X \).
\end{prop}

\begin{pf}
  Let  \( Z \subset X_{m}\) be  an \( {\Bbb{G}}_m \)-invariant 
  closed subset. 
  Then, we obtain:
  \[ \overline{\mu}_{Xm}(\bA^1\times Z)= Z .\] 
  On the other hand, 
  \( \overline{\mu}_{Xm}(\{0\}\times Z)=\sigma_{m}\circ \pi_{m}(Z) \)
  by \ref{invariant}.
  Therefore, as \( Z \) is closed, it follows that
   \[ Z \supset \overline{ \sigma_{m}\circ\pi_{m}(Z)}\supset 
   \sigma_{m}(\overline{\pi_{m}(Z)}),\] 
  which yields \( \pi_{m}(Z)\supset \overline{\pi_{m}(Z)} \).
\end{pf}

\begin{prop}
  Let \( f:X\to Y \) be a morphism of \( k \)-schemes of finite type.
  Then the canonical morphism \( f_{\infty}:X_{\infty}\to Y_{\infty} \) 
  is induced such that the 
  following diagram is commutative:
  \[ \begin{array}{ccc}
      X_{\infty}& \stackrel{f_{\infty}}\longrightarrow & Y_{\infty}\\
      \pi_{Xm} \downarrow\ \ \ \ \ & & \ \ \ \ \downarrow \pi_{Ym}\\
      X & \stackrel{f}\longrightarrow & Y\\
      \end{array}. \] 
\end{prop}

\begin{pf}
  The morphism \( f_{\infty} \) is induced as the projective limit of 
  \( f_{m} \)  \(( 
  m\in \bN) \).
\end{pf}

\begin{prop}
\label{proper}
%**(proper)
  Let \( f:X\to Y \) be a proper birational morphism of \( k 
  \)-schemes of finite type such that 
  \( f|_{X\setminus W}:X\setminus W\simeq Y\setminus V \), where \( W\subset X \) 
  and
  \( V\subset Y \) are closed.
  Then \( f_{\infty}  \) gives a bijection 
  \[ X_{\infty}\setminus W_{\infty} \to Y_{\infty}\setminus V_{\infty}. \]
\end{prop}

\begin{pf}
  Let \( \alpha \in Y_{\infty}\setminus V_{\infty} \), 
  then \( \alpha(\eta)\in X\setminus V \).
  As \( X\setminus W\simeq Y\setminus V \). 
  We obtain the following commutative diagram:
  \[ \begin{array}{ccc}
    \spec K((t))& \to & Y\\
    \downarrow & & \downarrow\\
    \sT & \stackrel{\alpha}\longrightarrow & X\\
    \end{array}. \]
    Then, as \( f \) is a proper morphism, 
    by the valuative criteria of properness, 
    there is a unique morphism \( \tilde \alpha: \sT \to Y \) such that 
    \( f\circ \tilde \alpha=\alpha \).
    This shows the bijectivity as required.
\end{pf}

  The following is the version for \( m=\infty \) of Proposition 
  \ref{etale}:
  
\begin{prop}
\label{etaleinfty}
%****(etaleinfty)
 If \( f:X\to Y \) is an \'etale morphism, then \( X_{\infty}\simeq 
 Y_{\infty}\times_{Y}X \).
\end{prop}

\begin{pf}
  As \( \varprojlim_{m}(Y_{m}\times_{Y}X)=(\varprojlim_{m}Y_{m})\times_{Y}X 
  \), the case \( m=\infty \) is reduced to the case \( m<\infty \) 
  which is proved in Proposition \ref{etale}.
\end{pf}

\begin{prop}
  There is a canonical isomorphism:\[ (X\times Y)_{m}\simeq 
  X_{m}\times  Y_{m}, \]
  for every \( m \in \bN\cup \{\infty\} \).
\end{prop}

\begin{pf}
  For an arbitrary \( k \)-scheme \( Z \),
  \[\Hom_{k}(Z, X_{m}\times Y_{m}) \simeq 
  \Hom_{k}(Z, X_{m})\times \Hom_{k}(Z, Y_{m}),
   \]
   and the right hand side is isomorphic to 
   \[ \Hom_{k}(Z\times_{\spec k}\stm, X)
   \times \Hom_{k}(Z\times_{\spec k}\stm, Y)\]
   \[\simeq \Hom_{k}(Z\times_{\spec k}\stm, X\times Y). \]
   \[ \simeq \Hom_{k}(Z, (X\times Y)_{m}). \]
  The case  \( m=\infty \) follows from this.
\end{pf}

\begin{prop}
  Let \( f:X\to Y \) be an open immersion (resp. closed immersion) of
  \( k \)-schemes of finite type.
  Then the induced morphism \( f_{m}:X_{m}\to Y_{m} \) is also
  an open immersion (resp. closed immersion) for every \( m\in 
  \bN\cup\{\infty\} \).
\end{prop}

\begin{pf}
  The open case  follows from Lemma \ref{open} and Proposition 
  \ref{etaleinfty}.
  For the closed case, we may assume that \( Y \) is affine.
  If \( Y \) is defined by \( f_{i} \) \( (i=1,.,r) \) in an affine 
  space, then \( X \) is defined by \( f_{i} \) \( (i=1,.,r,.,u) \)
  with \( r\leq u \) 
  in the same affine space.
  Then, \( Y_{m}  \) is defined by \( F^{(s)}_{i} \) \( (i=1,.,r,\  s\leq 
  m) \) and
  \( X_{m} \) is defined by  \( F^{(s)}_{i} \) \( (i=1,.,r,.,u,\  s\leq m) \) in 
  the corresponding affine space.
  This shows that \( X_{m} \) is a closed subscheme of \( Y_{m} \).
\end{pf}

\begin{rem}
  In the above proposition we see that the property open or closed 
  immersion of the base spaces  
   is inherited by the morphism of the space of jets and arcs.
  But some properties are not inherited. 
  For example, surjectivity and closedness are not inherited.  
\end{rem}

\begin{exmp}
  There is an example that \( f:X\to Y  \) is surjective and closed 
  but \( f_{\infty}:X_{\infty}\to Y_{\infty} \) is neither surjective 
  nor closed.
  Let \( X=\bA_{\bC}^2 \) and \( G=\langle \left( \begin{array}{cc}
                                   \epsilon & 0\\
                                    0 & \epsilon^{n-1}\\
                                   \end{array}\right)\rangle \) be a 
   finite cyclic subgroup in \( \GL(2, \bC) \) acting on \( X \),
   where \( n\geq 2 \) and \( \epsilon \) is a primitive \( n \)-th 
   root of unity.
  Let \( Y=X/G \) be the quotient of \( X \) by the action of \( G \).
  Then, it is well known that the singularity appeared in \( Y \) is
  \( A_{n-1} \)-singularity.    
  Then the canonical projection \( f:X\to Y \) is closed and 
  surjective. 
  We will see that these two properties are not inherited by 
  \( f_{\infty}:X_{\infty}\to Y_{\infty} \).   
  Let \( p \) be the image \( f(0)\in Y \).
  Then, by the commutativity 
  \[ \begin{array}{ccc}
  X_{\infty}& \stackrel{f_{\infty}}\longrightarrow& Y_{\infty}\\
  \ \ \ \downarrow \pi_{X}& &\ \ \ \downarrow \pi_{Y}\\
  X& \stackrel{f}\longrightarrow& Y,\\ 
  \end{array}\]
  we obtain \( \pi_{X}^{-1}(0)=f_{\infty}^{-1}\circ \pi_{Y}^{-1}(p) \).
  Here,  \( \pi_{X}^{-1}(0) \) is irreducible, since \( X \) is 
  non-singular.
  On the other hand \( \pi_{Y}^{-1}(p) \) has \( (n-1) \)-irreducible 
  components by \cite{nash}, \cite{i-k}.
  Therefore the morphism \( 
  f_{\infty} \)
  is not surjective for \( n\geq 3 \)¥.
  As \( X\setminus\{0\}\to Y\setminus \{p\} \) is \'etale, 
  The morphism 
  \[ (X\setminus\{0\})_{\infty}\to (Y\setminus \{p\})_{\infty} \] is 
  also \'etale by Proposition \ref{etaleinfty}.
  Since \( Y_{\infty} \) is irreducible,
  \( f_{\infty} \) is dominant.
  Therefore, \( f_{\infty} \) is not closed.
\end{exmp}

Next we think of the irreducibility of the arc space or jet 
schemes.
The following is proved in \cite{kln}.
In \cite{i0} we gave another proof by using \cite[Lemma 2.12]{i-k} and 
a resolution of the singularities.
Here we show a proof without a resolution.

\begin{thm}[\cite{kln}, \cite{i0}]
\label{good}
%**(good)
  If characteristic of \( k \) is zero, then the space of arcs of 
  a variety \( X  \) is irreducible.
\end{thm} 

\begin{pf}
  By \cite[Lemma 2.12]{i-k} we obtain the following: 
  \begin{enumerate}
  \item[(1)]
  Given any arc $\phi:\spec k'[[s]]\to X$, we construct an arc $\Phi$ such that $\phi\in \overline{\{\Phi\}}$
  and $\Phi(0)=\Phi(\eta)=\phi(\eta)$.
  \item[(2)]
  We construct an arc $\Psi$  such that  $\Phi\in \overline{\{\Psi\}}$ and 
  $\Psi(\eta)\in X\setminus \sing X$.
\end{enumerate}

Now for this $\Psi$ we apply the procedure (1) again, 
then we obtain a new arc $\Psi'$ such that  $\Psi\in \overline{\{\Psi'\}}$
  and $\Psi'(0)=\Psi'(\eta)=\Psi(\eta)\in X\setminus \sing X$.
  If we denote $\pi(\Psi')=\Psi'(0)=\lambda$, then $\Psi'\in 
  \pi^{-1}(\lambda) $.
  As $\lambda\in X\setminus \sing X$, it follows that 
  $$\Psi'\in \pi^{-1}(\lambda)\subset \overline{\pi^{-1}(\rho)},$$
  where \( \rho \) is the generic point of \( X \).
  This yields $\phi\in \overline{\pi^{-1}(\rho)}$ which is an 
  irreducible closed subset.
\end{pf} 

\begin{exmp}[\cite{i-k}, Example 2.13]
\label{charap}
  If the characteristic of \( k \) is \( p>0 \),
  \( X_{\infty} \) is not necessarily irreducible.
  For example, the hypersurface \( X \) defined by \( x^p-y^pz=0 \) 
  has an irreducible component in \( (\sing X)_{\infty} \) which is not in the 
  closure of \( X_{\infty}\setminus (\sing X)_{\infty} \).
\end{exmp}

\begin{exmp}[\cite{i}]
  Let \( X \) be a toric variety over an algebraically closed field of 
  arbitrary characteristic.
  Then, \( X_{\infty} \) is irreducible.
\end{exmp}

 Next let us think of \( m \)-jet scheme.
 A space of \( m \)-jets is not necessarily irreducible  even if the 
 characteristic of \( k \) is zero (see Example \ref{nonirred}).
 
 \begin{thm}[\cite{must01}] 
 If \( X \) is a variety of locally   complete intersection  over an 
 algebraically closed field of characteristic zero, 
 then \( X_{m} \) is irreducible for all \( m\geq 1  \)
 if and only if \( X  \) has rational singularities.
 \end{thm}

Another story in which a geometric property of space of jets 
determines the singularities on the base space is as follows:

\begin{thm}[\cite{ein}]
  Let \( X \) be a reduced divisor on a nonsingular variety over \( 
  \bC \).
  \( X \) has terminal singularities if and only if \( X_{m} \) is 
  normal for every \( m \in \bN\).  
\end{thm} 

%%%%%%%%%%%%%%%%%%%%%%%%%%%%%%%%%%%%%%%%%%%
%%%%%%%%%%%%%%%%%%%%%%%%%%%%%%%%%%%%%%%%%%%

\section{Introduction to the Nash problem}

\noindent
In this section, we assume the existence of resolutions of 
singularities. 
It is sufficient to assume that the characteristic of \( k \) is zero.
One of the most mysterious and fascinating problem in arc spaces is 
the Nash problem which was posed by Nash in his preprint in 1968.
It is a question about the Nash components and 
the essential divisors.
First we introduce the concept of essential divisors.

\begin{defn} Let $X$ be a variety, \( g:X_1 \to X \)
a proper birational morphism from a normal variety \( X_{1} \)
 and $E\subset X_1$ 
an irreducible exceptional divisor  of 
  \( g \). 
 Let  \( f:X_2 \to X \) be another
 proper birational morphism from a normal variety \( X_{2} \).
The birational  map
  \(f^{-1}\circ g  : X_1  \dasharrow  X_2 \) is 
defined on a (nonempty) open subset $E^0$ of $E$.
Because, by Zariski's main theorem, the ``undefined locus'' of 
a birational map between normal varieties 
is of codimension \( \geq 2 \).
The closure of $(f^{-1}\circ g)(E^0)$ is 
 called the {\it center} of $E$ on $X_2$.

  We say that \( E \) appears in \( f \) (or in \( X_2 \)),
  if 
the center of $E$ on $X_2$ is also a divisor. In this case
the birational  map
  \(f^{-1}\circ g  : X_1  \dasharrow  X_2 \) is  a local isomorphism at the 
  generic point of \( E \) and  we denote the birational
  transform of \( E \) on \( X_2 \)
  again by \( E \). For our purposes $E\subset X_1$ is identified
with $E\subset X_2$. 
 Such an equivalence class  is called an {\it exceptional divisor over $X$}. 
\end{defn}

\begin{defn}
  Let \( X \) be a variety over \( k \) and let \( \sing X \) be the 
  singular locus of \( X \). 
  In this paper, by  a
{\it  resolution} of the 
  singularities of \( X \) 
   we mean  a proper, birational  morphism \( f:Y\to X \) 
with  \( Y \)  non-singular
such that the restriction 
$Y\setminus f^{-1}(\sing X)\to X\setminus \sing X$ of \( f \)¥
 is an isomorphism.
\end{defn}

\begin{defn}\label{nashessential}
  An exceptional divisor \( E \) over \( X \) is called an {\it essential 
  divisor} over \( X \)
   if for every resolution \( f:Y\to X \)
the center of \( E \) on \( Y \) is an irreducible component of 
   \( f^{-1}(\sing X) \).

  For a given resolution \( f:Y\to X \), 
  the center of an essential divisor is called an {\it essential 
  component} on \( Y \).
\end{defn}  
  
\begin{prop}
\label{comp}
%*(comp)
  Let \( f:Y\to X \) be a resolution of the singularities of a 
  variety \( X \).  
  The set 
$$
{\cal E}={\cal E}_{Y/X}=
\left\{ 
    \begin{array}{c}
    \mbox{irreducible components of \( f^{-1}(\sing X)  \)}\\  
    \mbox{ which are  centers  of  essential divisors  over \( X \)}
    \end{array} 
    \right\}
$$
corresponds 
  bijectively to
  the set of all essential divisors  over \( X \).

  In particular, the set of essential divisors over \( X \) is a finite 
  set.
\end{prop}

\begin{pf}
  The map 
  \[\{\mbox{essential\ divisors\ over}\ X\}\to {\cal E}_{Y/X},\ \ 
  E\mapsto \mbox{center\ of\ } E \ \mbox{on}\ Y \]
  is surjective by the definition of essential components.
  To prove the injectivity, take  an essential component \( C \)
  and the  blow up \( Y'\to Y \) with the center
  \( C \).
  Then, there is a unique divisor \( E\subset Y' \) dominating \( C \).
  Let \( Y''\to Y' \) be a resolution of the singularities of \( Y' \).
  Then, \( E \) is the unique exceptional divisor on \( Y'' \) that 
  dominates  
  \( C \). 
  Therefore, every exceptional divisor over \( X \) with the center \( 
  C\subset Y \) has the center contained in \( E \) on a resolution 
  \( Y'' \) of the singularities of \( X \). 
  Therefore, by the definition of essential divisor, 
  this \( E \) is the unique essential divisor whose center 
  on \( Y \)
  is \( C \). 
\end{pf}

  C. Bourvier and G. Gonzalez-Sprinberg also introduce ``essential divisors'' 
  and ``essential components'' in \cite{B} and \cite{G-S},
  but we should note that the definitions are different from ours.
  In order to distinguish them we give different names to their 
  ``essential divisors'' and ``essential components''.
  
\begin{defn}[\cite{B}, \cite{G-S}]
  An exceptional divisor \( E \) over \( X \) is called a {\it BGS-essential 
  divisor} over \( X \)
  if \( E \) appears in every resolution.
  An exceptional divisor \( E \) over \( X \) is called a {\it BGS-essential 
  component} over \( X \) if the center of \( E  \) on every 
  resolution \( f \) of the singularity of \( X \) is an irreducible component 
  of \( f^{-1}(E') \), where \( E' \) is the center of \( E \) on \( X \).
\end{defn}  

We will see how different they are from our essential divisors and
essential components.
First we see that they coincide for 2-dimensional case. 
To show this we need to introduce the concept minimal resolution.

\begin{defn}
  A resolution \( f:Y\to X \) of the singularities of \( X \) is 
  called the {\it minimal resolution} 
  if for any resolution \( g:Y'\to X \), 
  there is a unique morphism \( Y'\to Y \) over \( X \).
\end{defn}
  
It is known that for a surface \( X \) the minimal resolution 
\( f:Y\to X \) exists. 
It is characterized that \( Y  \) has no exceptional curve of the 
first kind over \( X \).

For higher dimensional variety \( X \), the minimal resolution does 
not necessarily exist.
For example, \( X=\{xy-zw=0\}\subset \bA^4 \) has two resolutions 
neither of which dominates the other.
These two resolutions are obtained as follows: 
First take a blow-up \( f:\tilde Y \to X \) at the origin of \( X \) 
which has the unique 
singular point at the origin. 
Then, \( f \) is a resolution of the singularity of \( X \) and
 the exceptional divisor \( E \) of \( f \) is isomorphic to
  \( \bP^1\times \bP^1 \).
Here we have two contractions \( g_{1}:Y_{1}\to X \), \( g_{2}:Y_{2}
\to X \) whose restrictions are the first projection 
\( p_{1}:E=\bP^1\times \bP^1 \to \bP^1\)
and the second projection 
\( p_{2}:E=\bP^1\times \bP^1 \to \bP^1\), respectively.
The both \( Y_{i} \)'s are non-singular, therefore \( f_{i} \)'s are 
resolutions of the singularity of \( X \).
It is clear that there is no morphism between \( Y_{1} \) and 
\( Y_{2} \) over \( X \).

\begin{prop}
  If \( X \) is a surface, then each  set of
  ``essential divisors'',
  ``BGS-essential divisors'' and ``BGS-essential components'' are 
  bijective to 
  the set of  the components of the fiber \( f^{-1}(\sing X) \),
  where \( f:Y\to X \) is the minimal resolution.
  These are also essential components on the minimal resolution.
\end{prop}

\begin{rem}
  Four concepts ``essential divisor'', ``essential component'',
  ``BGS-essential divisor'' and ``BGS-essential component'' are 
  mutually different in general.
  
  First, our essential component is different from the others, 
  because it is a closed subset on a specific resolution and the others
  are all equivalence class of divisors.
 
  Next, a BGS-essential divisor is different from a BGS-essential 
  component or a essential divisor.
  Indeed, for \( X=(xy-zw=0)\subset \bA_{k}^4 \), the exceptional 
  divisor obtained by a blow-up at the origin is the unique essential 
  divisor and also the unique BGS-essential component, while 
  there is no BGS-essential divisor, since \( X \) has a  
  resolution whose exceptional set is \( \bP_{k}^1 \).
  
  Finally a BGS-essential component and an essential divisor are 
  different. Indeed, consider a cone generated by 
  \( (0,0,1),(2,0,1),(1,1,1), \) \((0,1,1) \) in \( \bR^3 \).
   It is well known that a cone generated by integer points in a real 
  Euclidean space defines an affine toric variety (see \cite{ful}, \cite{oda} for basic notion of 
  toric variety).
  Let \( X \) be the affine toric variety defined by this cone.
  Then the canonical subdivision adding a one dimensional cone 
  \( \bR_{\geq 0}(1,0,1) \) is a resolution of \( X \).
  As the singular locus of \( X \) is of dimension one, 
  there is no small resolution.
  Therefore, the divisor \( D_{(1,0,1)} \) is the unique 
  essential divisor, while \( D_{(1,1,2)} \) and \( D_{(2,1,2)} \) 
  are BGS-essential components by the criterion \cite[Theorem 2.3]{B}. 
\end{rem}

\begin{defn}
  Let \( X \) be a variety and \( \pi:X_{\infty}\to X \) the canonical 
  projection.
  An irreducible  component \( C \) of \( \pi^{-1}(\sing X) \) is 
  called a {\it Nash component} if it contains an arc \( \alpha  \) 
  such that \( \alpha(\eta)\not\in \sing X \).
  This is equivalent to that \( C \not\subset (\sing X)_{\infty} \).
\end{defn}

The following lemma is already quoted for the irreducibility of 
the space of arcs (Theorem \ref{good}).

\begin{lem}[\cite{i-k}]
   If the characteristic of the base field \( k \) is zero, 
   then every irreducible component of \( \pi^{-1}(\sing X) \) is 
   a Nash component. 
\end{lem}

  We note that for the positive characteristic case 
  this lemma does not hold. 
  Indeed, Example \ref{charap} is an example that 
   \( \pi^{-1}(\sing X) \)  has an irreducible component which is not
   a Nash component.

  Let \( f:Y\to X \) be a resolution of the singularities of \( X \)
  and \( E_{l} \) \( (l=1,..,r) \) the irreducible components of 
  \( f^{-1}(\sing X) \).
  Now we are going to introduce a map \( {\cal N} \) which is called 
  the Nash map
  $$
\left\{
\begin{array}{c}
\mbox{Nash components}\\
\mbox{of the space of arcs}\\
\mbox{of\  $ X$}
\end{array}
\right\}
 \stackrel{\cal N}\longrightarrow   
    \left\{ 
    \begin{array}{c}
    \mbox{essential}\\  
    \mbox{components}\\
    \mbox{ on $Y$}
    \end{array} 
    \right\}
\simeq
  \left\{ 
    \begin{array}{c}
    \mbox{essential}\\  
    \mbox{divisors}\\
    \mbox{ over $X$}
    \end{array} 
    \right\}.
$$

\begin{say}[construction of the Nash map]
  The  resolution \( f:Y\to X  \) induces a morphism
  \( f_{\infty}:Y_{\infty}\to X_{\infty}  \) of schemes. 
  Let \( \pi_{Y}:Y_{\infty} \to Y \) be the canonical projection.
  As \( Y \) is non-singular, \(  (\pi_{Y})^{-1}(E_{l})  \) is 
  irreducible for every \( l \).  
  Denote  by \(  (\pi_{Y})^{-1}(E_{l})^o  \) 
  the open subset of \( (\pi_{Y})^{-1}(E_{l}) \) consisting of the
   points  corresponding to   arcs \( \beta:\Spec 
   K[[t]]\to Y \) such that \( \beta(\eta)\not\in \bigcup_{l}E_{l} \).
  Let \( C_{i} \) \( (i\in I )\) be the Nash components of \( X \).
   Denote  by \( C_{i}^o  \) 
  the open subset of \( C_{i} \) consisting of the
   points  corresponding to   arcs \( \alpha:\Spec 
   K[[t]]\to X\) such that \( \alpha(\eta)\not\in \sing X \).
   As \( C_{i} \) is a Nash component, we have \( C_{i}^o \neq 
   \emptyset \).
  The  restriction of \( f_{\infty} \) gives  
  \[ f'_{\infty}:
  \bigcup_{l=1}^r (\pi_{Y})^{-1}(E_{l})^o \to \bigcup_{i\in I}C_{i}^o.\] 
By Proposition \ref{proper}, \( f'_{\infty} \) is surjective.
  Hence,        
 for each \( i\in I 
  \) there is  a unique \( l_{i} \) such that \( 1\leq l_{i}\leq r \)
  and the generic point \( 
  \beta_{l_{i}} \) of
  \( (\pi_{Y})^{-1}(E_{l_{i}})^o  \) is mapped to the generic point \( 
  \alpha_{i} \)
  of \( C_{i}^o \).
  By this correspondence \(  C_{i}\mapsto E_{l_{i}} \) we obtain a map
  \[ {\cal N}: \left\{
\begin{array}{c}
\mbox{Nash components}\\
\mbox{of the space of arcs}\\
\mbox{through $\sing X$}
\end{array}
\right\}
\longrightarrow   
    \left\{ 
    \begin{array}{c}
    \mbox{irreducible }\\  
    \mbox{components}\\
    \mbox{of  $f^{-1}(\sing X)$}
    \end{array} 
    \right\}.
\]
\end{say} 
 
\begin{lem} 
  The map \( {\cal N} \) is an injective map to the subset 
  \{ essential components on \( Y \)\}.
\end{lem} 

\begin{pf}
  Let \( {\cal N}(C_{i})= E_{l_{i}}  \).
  Denote the generic point of \( C_{i} \) by \( \alpha_{i} \) and
  the generic point of \( (\pi_{Y})^{-1}(E_{l}) \) by \( \beta_{l} \).
  If \( E_{l_{i}}=E_{l_{j}} \) for \( i\neq j \),
  then \( \alpha_{i}=f'_{\infty}(\beta_{l_{i}})= f'_{\infty}(\beta_{l_{j}})=
   \alpha_{j}\),
  a  contradiction.
  
  To prove that the \( \{E_{l_{i}}: i\in I\} \) are essential 
  components on \( Y \), let \( Y'\to X \) be another  resolution and
   \( \tilde Y \to X \) a divisorial resolution which
  factors through both \( Y \) and \( Y' \).
  Let \( E'_{l_{i}}\subset Y' \) and \( \tilde E_{l_{i}}\subset \tilde Y  \)
   be the irreducible components of the exceptional sets  
  corresponding to \( C_{i} \).
  Then, we can see that \(  E_{l_{i}}  \) and \(  E'_{l_{i}} \) are 
  the image of \( \tilde E_{l_{i}} \).
  This shows that \( \tilde  E_{l_{i}} \) is an essential divisor 
  over \( X \) and therefore \(  E_{l_{i}} \) is an essential 
  component on \( Y \).  
\end{pf}

\begin{problem}
   Is the Nash map 
    $$
\left\{
\begin{array}{c}
\mbox{Nash components}\\
\mbox{of the space of arcs}\\
\mbox{through $\sing X$}
\end{array}
\right\}
 \stackrel{\cal N}\longrightarrow   
    \left\{ 
    \begin{array}{c}
    \mbox{essential}\\  
    \mbox{components}\\
    \mbox{ on $Y$}
    \end{array} 
    \right\}
\simeq
  \left\{ 
    \begin{array}{c}
    \mbox{essential}\\  
    \mbox{divisors}\\
    \mbox{ over $X$}
    \end{array} 
    \right\}.
$$
 bijective?
\end{problem}
  After Nash's preprint  
   which posed this problem 
  was circulated in 1968, 
  Bouvier, Gonzalez-Sprinberg, 
  Hickel, Lejeune-Jalabert, Nobile, Reguera-Lopez and others 
  (see, \cite{B}, \cite{GL}, \cite{H}, \cite{L80}, \cite{L}, 
  \cite{LR}, \cite{nobile}, \cite{RL}) worked on the arc space of a 
  singular variety related to this problem.
 
  Recently for a  toric variety of arbitrary dimension  
  the Nash problem is affirmatively answered  
  but is negatively answered in general by Ishii and Koll\'ar in \cite{i-k}.

  Here, we show known results for this problem.
  
\begin{thm}[\cite{nash}]
  The Nash problem is affirmatively answered for  an \( A_{n} \)-singularity \( (n\in 
  \bN) \), where an \( A_{n} \)-singularity is the hypersurface singularity defined 
  by \( xy-z^{n+1}=0 \) in \( \bA_{k}^3 \).
\end{thm}  

\begin{thm}[\cite{RL}]
  The Nash problem is affirmatively answered for a minimal surface singularity.
  Here, a minimal surface singularity means a rational surface
  singularity with the reduced fundamental cycle.
  The fundamental cycle is induced by M. Artin (see  \cite{artin} for 
  the definition).
\end{thm}

\begin{thm}[\cite{LR}, \cite{rl03}]
\label{sand}
  The Nash problem is affirmatively answered for a sandwiched surface singularity.
  Here, a sandwiched surface singularity means the formal 
  neighborhood of a singular point on a surface obtained by blowing 
  up a complete ideal in the local ring of a closed point on a 
  non-singular algebraic surface. 
  A complete ideal is defined by O. Zariski and  Samuel (see 
  \cite{zs},  Vol II, Appendix 4 ), but the idea of a sandwiched 
  singularity is that it is a singularity which is birationally 
  sandwiched by non-singular surfaces. 
\end{thm}

These are results on rational surface singularities, 
the following gives affirmative answer for some non-rational surface  
singularities:

\begin{thm}[\cite{c-pp}]
  The Nash problem is affirmatively answered for a normal surface singularities with 
  the reduced fiber \( E \) of the singular point on the minimal resolution such that 
  \( E\cdot E_{i}<0 \) for every irreducible component \( E_{i}  \) of \( 
  E \).
\end{thm}

\noindent
This result is generalized to a wider class of surface singularities
in \cite{mo}.
We omit the statement, since it is not simple.

The following results are for arbitrary dimension.

\begin{thm}[\cite{i-k}]
  The Nash problem is affirmatively answered for a toric singularity of 
  arbitrary dimension.
\end{thm}

\begin{thm}[\cite{i2}]
  The Nash problem is affirmatively answered for a non-normal toric variety of 
  arbitrary dimension.
\end{thm}  

We have a notion of the local Nash problem which is a slight 
modification of the Nash problem (\cite{i3}).

\begin{thm}[\cite{i3}]
  The local Nash problem hold true for a quasi-ordinary singularities.
  Here, a quasi-ordinary singularity is a hypersurface singularity 
  which is a finite cover over a non-singular variety with the normal 
  crossing 
  branch locus. We note that a quasi-ordinary singularity is not 
  necessarily normal.
   
\end{thm}

The paper \cite{p-pp2} gives the affirmative answer to the Nash 
problem for a certain class of higher dimensional non-toric singularities.

So far we have seen the affirmative answers.
But there are negative examples  given in \cite{i-k}.

\begin{exmp}
  Let \( X \) be a hypersurface defined by \( x_{1}^3+
  x_{2}^3+x_{3}^3+x_{4}^3+x_{5}^6=0 \) in \( \bA_{\bC}^5 \).
  Then the number of the Nash components is one, while the 
  number of the essential divisors is two.
  Therefore the Nash map is not bijective.
\end{exmp}    

  By the above example we can construct  counter examples to the Nash problem  
   for any dimension greater than 3.
   At this moment the Nash problem is still open for two and three 
   dimensional variety.
   Now we can formulate a new version of the Nash problem:
   
\begin{problem}
  What is the image of the Nash map? 
  For two and three dimensional case, the image of the Nash map 
  coincides with the set of essential divisors?
\end{problem}   

 Related to this problem, we have one characterization of the image of 
 the Nash map given by  Reguera \cite{r}. 
 To formulate her result, we introduce the concept ``wedge''.
 
 \begin{defn}
 Let \( K\supset k \) be a  field extension. 
 A \( K \)-{\it wedge} of \( X \) is a \( k \)-morphism 
 \( \gamma:\spec K[[\lambda, t]]\to X \).
 A \( K \)-wedge \( \gamma \) can be identified to a \( K[[\lambda]] 
 \)-point on \( X_{\infty} \). 
 We call the {\it special arc} of \( \gamma \) the image in \( 
 X_{\infty} \) of the closed point \( 0 \) of \( \spec K[[\lambda]] \).
 We call the {\it generic arc} of \( \gamma \) the image in \( X_{\infty} 
 \) of the generic point \( \eta \) of \( \spec K[[\lambda]] \).
 
 \end{defn}.
 
\begin{thm}[\cite{r}]
  Let \( E \) be an essential divisor over \( X \) and \( f:Y\to X \) 
  a resolution of the singularities of \( X \) on which \( E \) appears.
  Let \( \alpha\in X_{\infty} \) be the generic point of \( 
  f_{\infty}(\pi_{Y}^{-1}(E)) \) and \( k_{E} \) the residue field of 
  \( \alpha \). 
  Then the following conditions are equivalent:
\begin{enumerate}
\item[(i)]
  \( E \) belongs to the image of the Nash map;
\item[(ii)]
  For any resolution of the singularities \( g:Y'\to X \) and for any 
  field extension \( K \) of \( k_{E} \), any \( K \)-wedge \( \gamma \) 
  on \( X \) whose special arc is \( \alpha \) and whose generic arc 
  belongs to \( \pi_{X}^{-1}(\sing X) \), lifts to \( Y' \);
\item[(iii)]
  There exists a resolution of the singularities \( g:Y'\to X \) 
  satisfying condition (ii).    
\end{enumerate}  
\end{thm}

As a corollary of this theorem, we also obtain 
Theorem \ref{sand}.

There are some notions ``the Nash problem for a pair \( (X, Z) \)'' 
consisting of a variety \( X \) and a closed subset \( Z \) (see 
\cite{petrov}, \cite{perez}). 
It seems that these are on the way of developing.

 %%%%%%%%%%%%%%%%%%%reference%%%%%%%%%%%%%%%%%%%%

\makeatletter \renewcommand{\@biblabel}[1]{\hfill#1.}\makeatother


\begin{thebibliography}{11}

 \bibitem{artin} M. Artin, {\em On isolated rational singularities of 
 surfaces}, Amer. J. Math. {\bf 88} (1966) 129--136.

\bibitem{B} C. Bouvier, {\em  Diviseurs essentiels, composantes essentielles
 des vari\'et\'es
toriques singuli\`eres},
Duke Math. J. {\bf 91} (1998) 609--620

\bibitem{G-S} C. Bouvier and G. Gonzalez-Sprinberg, {\em
Syst\`eme g\'en\'erateur minimal, diviseurs essentiels et 
G-d\'esingularisations de vari\'et\'es toriques}, Tohoku Math. J.
{\bf 47}, (1995) 125--149. 

\bibitem{blr} S. Bosch, W. L\"utkebohmert and M. Raynaud, 
{\em N\'eron Models}, Ergebnisse der Mathematik und ihrer 
Grenzgebiete, {\bf 21} (1990) Springer-Verlag.



\bibitem{craw} A. Craw, {\em An introduction to motivic integration,}
math.AG/9911179 

\bibitem{DEI} T. De Fernex, L. Ein and S. Ishii, {\em Divisorial 
valuations via arcs}, preprint (2007)  math.AG/0701867.


\bibitem{DL1} J. Denef and F. Loeser, {\em Motivic Igusa 
zeta-functions,} J. Alg. Geom. {\bf 7} (1988) 505--537.



\bibitem{DL2} J. Denef and F. Loeser, {\em Germs of arcs on singular
varieties and motivic integration,} Invent. Math. {\bf 135}, (1999)
201--232. 

\bibitem{DL3} J. Denef and F. Loeser, {\em Motivic exponential 
integrals and a motivic Thom-Sebastiani Theorem,} Duke Math. J. {\bf 
99} (1999) 201--232.

\bibitem{DL4} J. Denef and F. Loeser, {\em Motivic integration, 
quotient singularities and the McKay correspondence,} Comositio Math.
{\bf 131}, (2002) 267-290.
 
\bibitem{DL5} J. Denef and F. Loeser, {\em Motivic integration and the 
Grothendieck group of pseudo-finite fields,} Proceeding of the 
International Congress of Mathematicians (Beijing, 2002) vol II, 
Higher Ed. Press, Beijing, 2002, 13--23.



\bibitem{ein} L. Ein, M. Musta\cedilla{t}\v{a} and T. Yasuda,
{\em Jet schemes, log discrepancies and inversion of adjunction,}
 Invent. Math. {\bf 153} (2003) 519-535.

\bibitem{e-Mus} L. Ein and M. Musta\cedilla{t}\v{a}. {\em Inversion of 
Adjunction for local complete intersection varieties}, Amer. J. 
Math. 126 (2004), 1355--1365.

\bibitem{ELM} L. Ein, R. Lazarsfeld and M. Musta\cedilla{t}\v{a},
{\em Contact loci in arc spaces}, Compositio Math. {\bf 140} (2004) 
1229--1244. 

\bibitem{ful} W. Fulton, {\em Introduction to Toric Varieties}, Annals 
of Math. Studies {\bf 131}, Princeton Univ. Press. (1993).

\bibitem{perez} P.D. Gonz\'alez P\'eres, {\em Nash problem for quasi-ordinary hypersurface 
singularities} preprint 2006.

\bibitem{GL} G. Gonzalez-Sprinberg and M. Lejeune-Jalabert, {\em Families of smooth curves on
surface singularities and wedges}, Annales Polonici Mathematici, 
LXVII.2, (1997) 179--190.

\bibitem{gr} M. Greenberg, {\em Rational points in henselian discrete 
valuation rings,} Publ. Math. I.H.E.S. {\bf 31} (1966), 59--64.

\bibitem{ha} R. Hartshorne, {\em Algebraic Geormetry}, 
Graduate Texts in Math. {\bf 52} Springer-Verlag, (1977).

\bibitem{H} M. Hickel, {\em Fonction de Artin et germes de courbes
trac\'ees sur un germe d'espace analytique}, Amer. J. Math. {\bf 115},
(1993) 1299--1334.

\bibitem{i-k} S. Ishii and J. Koll\'ar, {\em The Nash problem   on  
arc families of singularities,}  Duke 
Math. J. {\bf 120} No.3 (2003) 601-620.

\bibitem{i0} S. Ishii, {\em Introduction of arc spaces and the Nash 
problem,} RIMS Kokyu-Roku, {\bf 1374}, (2004) 40--51.


\bibitem{i} S. Ishii, {\em The arc space of a toric variety}, J. 
Algebra,
{\bf 278} (2004) 666--683

\bibitem{i2} S. Ishii, {\em Arcs, valuations and the Nash map}, J. 
reine angew. Math, {\bf 588} (2005) 71--92.

\bibitem{i3} S. Ishii, {\em The local Nash problem on arc families of 
singularities},  Ann. Inst. Fourier, Grenoble {\bf 56} (2006) 1207-1224.

\bibitem{i4} S. Ishii, {\em Maximal divisorial sets in arc spaces}, 
to appear in Proceeding of Alg. Geom. in East Asia II 2005, Advanced 
Studies in Pure Math..



\bibitem{kln} E. R.  Kolchin,
{\em Differential algebra and algebraic groups},
 Pure and Applied Mathematics, Vol. {\bf 54}, Academic Press, New York-London, 1973.



\bibitem{ko} M. Kontsevich, {\em Lecture at Orsay} (December 7, 1995)

\bibitem{L80} M. Lejeune-Jalabert, {\em Arcs analytiques et
r\'esolution minimale des 
surfaces quasihomog\`enes}. in: Lecture Notes in  Math. {\bf 777}, (1980) 
303--336.

\bibitem{L} M. Lejeune-Jalabert, {\em Courbes trac\'ees sur un germe 
d'hypersurface}. Amer. J. Math. {\bf 112}, (1990) 525--568.

\bibitem{LR} M. Lejeune-Jalabert and A. J. Reguera-Lopez, {\em Arcs 
and wedges on sandwiched surface singularities}, Amer. J. Math. 
{\bf 121}, (1999) 1191--1213. 

\bibitem{l} F. Loeser, {\em Seattle lectures on motivic integration}. preprint 2006 in his web page.

\bibitem{mo} M. Morales, {\em The Nash problem on arcs for surface 
singularities}. preprint (2006) math.AG/0609629

\bibitem{must01} M. Musta\cedilla{t}\v{a}, {\em Jet schemes of locally 
complete intersection canonical singularities,} with an appendix by 
David Eisenbud and Edward Frenkel, Invent. Math. {\bf 145} (2001) 
397--424.

\bibitem{must02} M. Musta\cedilla{t}\v{a}, {\em
Singularities of Pairs via Jet Schemes},
 J. Amer. Math. Soc. {\bf 15} (2002), 599--615.

\bibitem{nash} J. F. Nash, {\em Arc structure of singularities},
Duke Math. J. {\bf 81}, (1995) 31--38. 

\bibitem{nobile} A. Nobile, {\em On Nash theory of arc structure of singularities}, Ann. Mat. Pura
Appl. {\bf 160} (1991), 129--146.

\bibitem{oda} T. Oda, {\em Convex Bodies and Algebraic Geometry}, 
Ergeb. Math. Grenzgeb. {\bf 15}, Springer-Verlag  (1988).

\bibitem{petrov} P. Petrov, {\em Nash problem for stable toric 
varieties,} preprint (2006) math.AG/0604432.


\bibitem{c-pp} C. Pl\'enat and P. Popescu Pampu, {\em A class of 
non-rational surface singularities for which the Nash map is bijective}
math.AG/0410145.


\bibitem{p-pp2} C. Pl\'enat and P. Popescu-Pampu,
{\em Families of higher dimensional germs with bijective Nash map}, math.AG/0605566 


\bibitem{RL} A. J. Reguera-Lopez, {\em Families of arcs 
on rational surface singularities}, Manuscr. Math. 
{\bf 88}, (1995) 321--333.

\bibitem{rl03} A. J. Reguera-Lopez, {\em Image of the Nash map in 
terms of wedge,} C.R. Acad. Sci. Paris, Ser. I,  {\bf 338} (2004) 
385--390.

\bibitem{r} A. J. Reguera, {\em A curve selection lemma in spaces of 
arcs and the image of the Nash map}, Compositio Math. {\bf 142} (2006)
119--130.

\bibitem{v1} W. Veys, {\em Zeta functions and `Kontsevich invariants' on 
singular varieties,} Canadian J. Math. {\bf 53} (2001) 834--865.

\bibitem{v2} W. Veys, {\em Stringy invariants of normal surfaces,} J. 
Alg. Geom. {\bf 13} (2004) 115--141.

\bibitem{v3} W. Veys, {\em Stringy zeta functions of \( \bQ 
\)-Gorenstein varieties,} Duke Math. J. {\bf 120} (2003) 469--514. 


\bibitem{veys} W. Veys, {\em Arc spaces, motivic integration and 
stringy invariants}, math.AG/0401374. 

\bibitem{voj} P. Vojta, {\em Jets via Hasse-Schmidt derivations}, 
math.AG/0407113.

\bibitem{zs} O. Zariski and P. Samuel, {\em Commutative Algebra I, 
 II},
 Van Nostrand (1958), (1970)

\end{thebibliography}
\end{document}